\documentclass[12pt]{amsart}
\usepackage{amsmath}
\usepackage{amsfonts}
\usepackage{amssymb}

\newcommand{\ZZ}{\mathbb{Z}}
\newcommand{\FF}{\mathbb{F}}
\newcommand{\CC}{\mathbb{C}}
\newcommand{\QQ}{\mathbb{Q}}
\newcommand{\OO}{\mathcal{O}}
\newcommand{\sE}{\mathcal{E}}
\newcommand{\sT}{\mathcal{T}}
\newcommand{\A}{\mathbf{A}}
\newcommand{\NN}{\mathbf{N}}

\renewcommand{\bar}{\overline}
\newcommand{\Ehat}{\widehat{E}}

\newcommand{\hhat}{\hat{h}}
\newcommand{\Hhat}{\widehat{H}}

\DeclareMathOperator{\Aut}{Aut}
\DeclareMathOperator{\den}{den}
\DeclareMathOperator{\Gal}{Gal}
\newcommand{\NT}{\mathrm{NT}}
\DeclareMathOperator{\ord}{ord}
\newcommand{\sep}{\mathrm{sep}}

\newcommand{\infbar}{{\overline{\infty}}}
\newcommand{\laurent}[2]{#1(\!( #2 )\!)}
\newcommand{\mult}[1]{[#1]}
\newcommand{\pair}[2]{\langle #1,#2 \rangle}

\newcommand{\pairing}{\langle\phantom{x},\phantom{x}\!\rangle}

\newcommand{\pos}[1]{\left\langle #1 \right\rangle}
\newcommand{\power}[2]{#1[\![ #2 ]\!]}
\newcommand{\rhoKv}{\rho^K_v}
\newcommand{\rhokv}{\rho^k_v}
\newcommand{\smfrac}[2]{{\textstyle \frac{#1}{#2}}}
\newcommand{\sbstar}[2]{#1_{#2}^{\times}}

\newcounter{alphenum}
\newenvironment{alphenumerate}{\setcounter{alphenum}{0} \begin{list}
        {\textnormal{(\alph{alphenum})}}{\usecounter{alphenum}
                \setlength{\leftmargin}{1truecm}
                \setlength{\labelwidth}{7truemm}
                \setlength{\labelsep}{2truemm}
                \setlength{\itemsep}{2pt}
                \setlength{\parsep}{0pt}
                \setlength{\itemindent}{0pt}
        }}{\end{list}}

\makeatletter\@addtoreset{equation}{section}\makeatother

\newtheorem{thm}[equation]{Theorem}
\newtheorem{lem}[equation]{Lemma}
\newtheorem{cor}[equation]{Corollary}
\newtheorem{prop}[equation]{Proposition}
\newtheorem{hyp}[equation]{Hypothesis}
\theoremstyle{definition}
\newtheorem{rmk}[equation]{Remark}

\begin{document}

%%%%%%%%%%%%%%%%%%%%%%%%%%%%%%%%%%%%%%%%%
%               Opening                 %
%%%%%%%%%%%%%%%%%%%%%%%%%%%%%%%%%%%%%%%%%

  \title[Canonical heights in characteristic~$p$]{Canonical
        heights on elliptic curves in characteristic~$p$}
  \author{Matthew A.~Papanikolas}
  \address{Department of Mathematics\\
        Pennsylvania State University\\
        University Park, PA 16802}
  \email{map@math.psu.edu}

\begin{abstract}
Let $k = \FF_q(t)$ be the rational function field with finite constant
field and characteristic~$p \geq 3$, and let $K/k$ be a finite
separable extension.  For a fixed place~$v$ of~$k$ and an elliptic
curve~$E/K$ which has ordinary reduction at all places of~$K$
extending~$v$, we consider a canonical height pairing $\pairing_v
\colon E(K^{\sep}) \times E(K^{\sep}) \to \sbstar{\CC}{v}$ which is
symmetric, bilinear and Galois equivariant.  The pairing
$\pairing_\infty$ for the ``infinite'' place of~$k$ is a natural
extension of the classical N\'eron-Tate height.  For $v$ finite, the
pairing $\pairing_v$ plays the role of global analytic $p$-adic
heights.  We further determine some hypotheses for the non-degeneracy
of these pairings.
\end{abstract}

%%%%%%%%%%%%%%%%%%%%%%%%%%%%%%%%%%%%%%%%%
%               Document Text           %
%%%%%%%%%%%%%%%%%%%%%%%%%%%%%%%%%%%%%%%%%

\maketitle

\section{Introduction} \label{sec:intro}

Our goal in this paper is to construct and investigate canonical
heights on elliptic curves defined over global function fields in
characteristic~$p$.  These heights will take values in the completions
of the function field and serve as analogues of both the
classical N\'eron-Tate height and also analytic $p$-adic heights on
elliptic curves over number fields,
cf.~\cite{bern81,maztat83,MTT86,pr82,schn82}.

Except in Sections~\ref{sec:sigma} and~\ref{sec:sigform}, we maintain
the following notation throughout the paper.  We take~$p$ to be an odd
prime; $\FF_q$ the finite field with $q = p^n$ elements; $k$ the
rational function field~$\FF_q(t)$; and~$A$ the polynomial
ring~$\FF_q[t]$.  The valuations of the field~$k$ corresponding to
maximal ideals of~$A$ are called \emph{finite places}; the unique
remaining place is the \emph{infinite place,} denoted~$\infty$, with
$\ord_\infty = -\deg$, where $\deg$ is the degree map on polynomials
extended to rational functions.

We let $K/k$ denote a \emph{global function field,} i.e., a finite
separable field extension of~$k$, and we let~$\OO$ be the integral
closure of~$A$ in~$K$.  For a place~$w$ of~$K$, we let~$K_w$ denote
the completion at~$w$; $\OO_w$ its valuation ring; $\FF_w$ its residue
field; and~$\CC_w$ the completion of an algebraic closure of~$K_w$.
Furthermore, $U^1(K_w)$ and $U^1(\CC_w)$ are the groups of $1$-units.
A place of~$K$ is either finite or infinite depending on the place it
extends from~$k$.  Finally, for any field~$F$, we let $F^{\sep}$ be a
separable closure and~$\,\bar{\! F}$ an algebraic closure.

For a fixed place~$v$ of~$k$, we consider an elliptic curve~$E/K$
which has \emph{ordinary reduction} at all places of~$K$
extending~$v$.  That is, the formal group of each reduced curve has
height~$1$.  In Section~\ref{sec:Can} we define a Galois equivariant
quadratic form
$$
        \Hhat_v\colon E(K^{\sep}) \to \sbstar{\CC}{v},
$$
which equivalently induces a symmetric bilinear pairing,
$$
        \pairing_v\colon E(K^{\sep}) \times E(K^{\sep})
                \to \sbstar{\CC}{v}.
$$
We construct~$\Hhat_v$ as a product of local factors.  At the places
not dividing~$v$, these local factors are derived from intersection
multiplicities on the N\'eron model for~$E/K$, and at the places
above~$v$, the local heights are defined using the Mazur-Tate sigma
function~\cite{maztat91}.

It should be noted that these constructions contain differences
depending on whe\-ther $v$ is a finite or infinite place of~$k$.  We
find in Section~\ref{sec:NText} that for $v = \infty$, the resulting
$\infty$-adic height extends the N\'eron-Tate height.  Indeed we find
for all $P \in E(K^{\sep})$ that
$$
        \deg (\Hhat_\infty (P)) = 2\hhat_{\NT}(P).
$$
If $v$ is a finite place, the $v$-adic height takes values
in~$U^1(\CC_v)$.  In Sections~\ref{sec:nondeg} and~\ref{sec:kummer} we
investigate the non-degeneracy of these height functions.  Just as for
$p$-adic heights for elliptic curves defined over number fields, the
non-degeneracy of the~$v$-adic height is related to the non-existence
of universal norms coming from the so-called Carlitz cyclotomic tower
of~$K$ (see \cite{map:Univ,map:PhD}).

The assumption that $p$ is odd is not necessary, but for the purposes
of length we have disregarded the case where $p=2$.  Moreover, the
flavor of the results are identical~\cite{map:PhD}.

Finally, it is important to remark that the canonical heights
discussed in this paper are effectively computable, and in
Section~\ref{sec:exmp} we discuss an explicit example.

\section{Elliptic curves over global function fields} \label{sec:NT}

Fixing a field~$K/k$ which is the function field of a smooth curve
$X/\FF_q$, we let $E/K$ be an elliptic curve defined over~$K$ and
take~$\sE/X$ to be a N\'eron model for~$E/K$ with identity
component~$\sE^0/X$.  Using the intersection multiplicities of
sections on~$\sE$, it is possible to construct the classical
\emph{N\'eron-Tate canonical height pairing,}
$$
        \pairing_{\NT} \colon E(K^{\sep}) \times E(K^{\sep}) \to \QQ,
$$
which is symmetric and bilinear and when restricted to~$E(K)$ is
non-degenerate modulo torsion~\cite{sil:ATAEC}.  Equivalently, we can
construct the associated quadratic form~$\hhat_{\NT}$, which we define
so that $\hhat_{\NT}(P) = \frac{1}{2}\pair{P}{P}_{\NT}$.  We note that
here the N\'eron-Tate height is normalized to take values which are
independent of the chosen field of definition.

As is well known, the N\'eron-Tate height can be computed as a sum of
local heights (see~\cite{lang:ECDA,sil:ATAEC}).  We fix a Weierstrass
equation for~$E$,
\begin{equation} \label{eqn:wefix}
        y^2 + a_1xy + a_3y = x^3 + a_2x^2 + a_4x + a_6,\quad a_i \in
        K.
\end{equation}
For a place~$w$ of~$K$, we let $E_0(K_w) \subset E(K_w)$ denote the
subgroup of points whose reduction is non-singular at~$w$, i.e., the
points whose sections meet~$\sE^0$ on the fiber above~$w$.  There is a
local height function, $\lambda_w \colon E_0(K_w) \to \ZZ$, defined so
that for $P \in \sE^0(X) \subset E(K)$, the N\'eron-Tate height is
obtained by the formula
\begin{equation} \label{eqn:NTloc}
        \hhat_{\NT}(P) = \frac{1}{2\,[K:k]} \sum_{w \in X} [\FF_w :
        \FF_q] \lambda_w(P).
\end{equation}
If $P \in \sE^0(X)$ and (\ref{eqn:wefix}) is minimal at~$w$, then
$\lambda_w(P) = \max \{ -\ord_w(x(P)), 0 \}$.

\section{The Mazur-Tate sigma function} \label{sec:sigma}

Here we review some facts about the Mazur-Tate sigma function as
defined in~\cite{maztat91}.  This function, defined on the formal
group of an elliptic curve over a local field, is a natural
non-archimedean analogue of the classical Weierstrass sigma function,
and similar to the classical sigma function, it provides an analytic
local height function on an elliptic curve.  For more details on its
connection with $p$-adic heights for elliptic curves over number
fields see~\cite{MTT86}. There is also an analogue of the Weierstrass
zeta function, defined by Voloch~\cite{voloch:zeta}, which is also of
interest here.

Throughout this section and the next we diverge from the previous
notation and maintain the following conventions: $K$ is a field
complete with respect to a discrete valuation~$\ord$; $\OO$ is the
valuation ring of~$K$; and~$\FF$ is its residue field.  We assume that
both~$\OO$ and~$\FF$ have characteristic~$p$ and begin with some
necessary definitions.

Let $E/K$ be an elliptic curve with ordinary reduction, and fix an
invariant differential~$\omega$ on~$E$.  We choose a minimal
Weierstrass equation~(\ref{eqn:wefix}) so that $\omega = dx/(2y + a_1x
+ a_3)$.  Taking~$\Ehat/\OO$ for the formal completion of the N\'eron
model of~$E$ along its zero-section, we then pick a uniformizing
parameter~$z$ on~$\Ehat$ and let $\beta = (\omega/dz)(O)$.  The
multiplication-by-$p$ map on $\Ehat$ then has the form
\begin{equation} \label{eqn:Xp}
        \mult{p}(z) = \alpha\beta^{p-1}z^p + \cdots, \quad \in
        \power{\OO}{z^p},
\end{equation}
where $\alpha$ is the \emph{Hasse invariant} of $E/K$ (see~{\S}12.4
of~\cite{katmaz}).

The fraction field~$L$ of the power series ring~$\power{\OO}{z}$ is
naturally considered the field of rational functions on~$\Ehat$, and
it contains~$K(E)$ as a subfield.  For an integer $m$, the
\emph{$m$-th division polynomial with respect to~$\omega$} is the
function $f_m \in K(E) \subset L$ with divisor $\mult{m}^{-1}(O) -
m^2(O)$ such that
$$
        \frac{z^{m^2}\, f_m(z)}{\mult{m}(z)}(O) = \beta^{1-m^2},
$$
as defined in~\cite{maztat91}.  Such division polynomials satisfy many
recursion formulas and possess a rich structure.  Notably,
\begin{equation} \label{eqn:fpz}
        f_{p^n}(z) = f_{p^{n-1}}( \mult{p}(z)) f_p(z)^{p^{2n-2}}.
\end{equation}

\begin{thm}[Mazur-Tate~\cite{maztat91}] \label{thm:MT}
The sigma function $\sigma = \sigma_{E, \omega}$ is characterized as
the unique function in $\beta z(1 + z\power{\OO}{z})$ satisfying any
one of the following equivalent properties.
\begin{alphenumerate}
\item If $P$, $Q$ are nonzero points in $\Ehat(\OO)$, then
$$
        \frac{\sigma(P-Q)\sigma(P+Q)}{\sigma(P)^2 \sigma(Q)^2}
                = x(Q) - x(P).
$$
\item If $m \in \ZZ$ and $Q \in \Ehat(\bar{\OO})$, then
$$
        \sigma(mQ) = \sigma(Q)^{m^2}\, f_m(Q),
$$
where $f_m$ is the $m$-th division polynomial with respect to~$\omega$.
\item If $Q \in \Ehat(\OO)$, then
$$
        \sigma(pQ) = \sigma(Q)^{p^2}\, f_p(Q).
$$
\end{alphenumerate}
\end{thm}

\section{A formula for the sigma function} \label{sec:sigform}

Recalling~(\ref{eqn:Xp}), we deduce that $\mult{p^n}(z) = \alpha_n
\beta^{p^n-1} z^{p^n} + \cdots \in\power{\OO}{z^{p^n}}$, where
$\alpha_n = \alpha \cdot \alpha^p \cdot \alpha^{p^2} \cdot\,\cdots\,
\cdot \alpha^{p^{n-1}} = \alpha^{(p^n-1)/(p-1)}$.  An exercise in
formal power series shows that inverting~$\mult{p^n}(z)$ with respect
to composition yields
\begin{equation} \label{eqn:Xpinv}
        \mult{p^n}^{-1}(z) = \frac{1}{\beta} \left(
                \frac{\beta}{\alpha_n} \right)^{1/p^n} z^{1/p^n}
                + \cdots \quad \in \power{\OO^{1/p^n}}{z^{1/p^n}},
\end{equation}
which when raised to the $p^n$-th power results in a series
in~$\power{\OO}{z}$.  This series will converge for points in the
formal group, which highlights the fact that the points in $\Ehat$ are
uniquely $p$-divisible since $\Ehat$ has height~$1$.

\begin{thm} \label{thm:sigprod}
Let $\sigma(z)$ be the power series given by the infinite product
$$
        \sigma(z) = \beta z \, \prod_{n \geq 1} \left[
                \frac{\beta^{p-1}}{\alpha^{p^{n-1}}\,
                \alpha_n^{p-1}} \, z^{p-1} f_p (
                \mult{p^n}^{-1}(z) )^{p^{n-1}} \right]^{p^{n-1}}.
$$
The product converges in $\power{\OO}{z}$ and the resulting series is
the Mazur-Tate sigma function~$\sigma_{E,\omega}$.
\end{thm}

\begin{proof}
For convergence, it suffices to show that each factor in the square
brackets above is a $1$-unit in $\power{\OO}{z}$, which follows from
(\ref{eqn:Xpinv}).  A straightforward manipulation of the factors
permits us to prove that $\sigma(z)$ satisfies the identity
in~(\ref{thm:MT}c) and thus represents the sigma function.
See~\cite{map:PhD} for further details.
\end{proof}

\section{Canonical heights} \label{sec:Can}

We return to the notation given in Sections~\ref{sec:intro}
and~\ref{sec:NT}.  For a global function field $K/k$, let $E/K$ be an
elliptic curve with ordinary reduction at all places of~$K$ which
extend a chosen place~$v$ of~$k$.  In this section we construct a
symmetric bilinear pairing
\begin{equation} \label{eqn:par}
        \pairing_v \colon E(K) \times E(K) \to
        \sbstar{\CC}{v},
\end{equation}
which will serve as a canonical height.  Our interests will lie with
investigating various properties of these heights as well with
determining their non-degeneracy (see Sections~\ref{sec:nondeg}
and~\ref{sec:kummer}).  Although the constructions are similar, the
resulting pairing will have a different meaning depending on whether
the preferred place~$v$ is finite or infinite.  The $\infty$-adic
height plays the role of the N\'eron-Tate height (see
Theorem~\ref{thm:NText}) and the $v$-adic height ($v$ finite) is the
analogue of the analytic $p$-adic heights for elliptic curves over
number fields, cf.~\cite{maztat83,MTT86,schn82}.

The techniques in this section are fairly standard, so we merely
address the important points.  Moreover, our construction is similar
to that of~\cite{MTT86}; for more details the reader is directed
to~\cite{map:PhD}.  For each place $v$ of $k$ we fix a uniformizer
$\pi_v$ in $\FF_v k$, and we observe that the group $\pi_v^{\QQ}\cdot
U^1(\CC_v)$ is uniquely divisible.  For $x \in \sbstar{\CC}{v}$, we
let $\pos{x} = \pos{x}_v$ denote the \emph{positive part} of~$x$,
i.e., the class of $x$ in $\pi_v^{\QQ}\cdot U^1(\CC_v)$
modulo~$\sbstar{\bar{\FF}}{v}$.  The values of our heights will be
taken in the positive elements of $\sbstar{\CC}{v}$.

We then fix a Weierstrass equation for~$E/K$ as in~(\ref{eqn:wefix})
and let $\omega = dx/(2y + a_1x + a_3)$.  For each place~$w$ of~$K$,
we let $z_w$ be a local uniformizing parameter for $\Ehat/\OO_w$ and
$\beta_w = (\omega/dz_w)(O) \in \sbstar{K}{w}$.  We consider the
subgroup $E_v(K) \subset E(K)$ of finite index consisting of points on
the identity component~$\sE^0$ of~$\sE$ which also specialize to~$O$
on the fibers of places which extend~$v$, i.e.,
$$
        E_v(K) = \sE^0(X) \cap \bigcap_{w\mid v} \Ehat(\OO_w).
$$
For each $P \in E_v(K) \setminus \{ O\}$, we define an idele $i(P) \in
\sbstar{\A}{K}$ component-wise as follows:
\begin{equation} \label{eqn:I}
        i(P)_w = \left\{
                \begin{array}{ll}
                1 & \text{if $w \mid \infty$ and $v \neq \infty$,} \\
                \beta_w & \text{if $w \nmid v$ finite and $P \notin
                        \Ehat(\OO_w)$,} \\
                \beta_w z_w(P) & \text{if $w \nmid v$ finite and $
                        P \in \Ehat(\OO_w)$,} \\
                \sigma_{E/\!K_w,\omega}(P) & \text{if $w \mid v$.}
                \end{array}
        \right.
\end{equation}
We set $i(O) = 1$.  We observe that if $v$ is finite then $i(P)$ is
supported away from~$\infty$.  Furthermore, we find that for $w \nmid
v$ or~$\infty$,
\begin{equation} \label{eqn:ordiP}
        \ord_w(i(P)_w) = \smfrac{1}{2} \lambda_w(P).
\end{equation}
We then define a map $\rhoKv\colon \sbstar{\A}{K} \to \sbstar{\CC}{v}$
as follows.  We take
\begin{equation} \label{eqn:rhoKv}
        \rhoKv ((e_w)_w) = \rhokv \bigl( \NN^K_k (e_w)_w
                \bigr)^{1/[K:k]},
\end{equation}
where $\NN^K_k$ is the norm map on ideles and
$$
        \rhokv ((e_{\tilde{v}})_{\tilde{v}}) = \frac{1}{\pos{e_v}}
                \prod_{\tilde{v} \neq \infty}
                \pi_{\tilde{v}}^{\ord_{\tilde{v}}(e_{\tilde{v}})},
$$
for $\pi_{\tilde{v}}$ the positive uniformizer of $\tilde{v}$ in
$\FF_v A$.  Note that the values of $\rho^k_\infty$ are in
$\pi_{\infty}^{\QQ} \cdot U^1(\CC_{\infty})$, whereas for $v$ finite,
$\rhokv$ takes values in~$U^1(\CC_v)$.

By combining~(\ref{eqn:I}) and~(\ref{eqn:rhoKv}), we define a height
function by
$$
        \Hhat_v = \rhoKv \circ i^2 \colon E_v(K) \to \sbstar{\CC}{v},
$$
which is a quadratic form and hence a ``canonical height.''  We
observe that the values of $\Hhat_v$ are independent of the choices of
formal parameters and differentials made in~(\ref{eqn:I}).  Because
the image is contained in a uniquely divisible subgroup
of~$\sbstar{\CC}{v}$, we can extend the definition to all of~$E(K)$.
The following proposition summarizes the above discussion.

\begin{prop}
Let $K/k$ be a global function field and~$E/K$ an elliptic curve with
ordinary reduction at all the places of~$K$ extending~$v$.  There is a
unique quadratic form $\Hhat_v\colon E(K) \to \sbstar{\CC}{v}$ which
extends
$$
        \Hhat_v = \rhoKv\circ i^2 \colon E_v(K) \to \sbstar{\CC}{v}.
$$
Furthermore, the height behaves well under base-extension, and gives
rise to a symmetric, bilinear pairing
$$
        \pairing_v \colon E(K^{\sep}) \times E(K^{\sep}) \to
                \sbstar{\CC}{v},
$$
with $\Hhat_v(P) = \pair{P}{P}_v$, which is Galois equivariant, i.e.,
$\pair{\tau P}{\tau Q}_v = \pair{P}{Q}_v$ for all $P$,~$Q \in
E(K^{\sep})$ and all $\tau \in \Gal(K^{\sep}\!/K)$.
\end{prop}

\section{Extension of the N\'eron-Tate height} \label{sec:NText}

In this section, we will prove the following theorem which shows that
the \mbox{$\infty$-adic} height as defined in the previous section
actually extends the N\'eron-Tate height in a natural way.

\begin{thm} \label{thm:NText}
Let $E/K$ be an elliptic curve with ordinary reduction at all the
infinite places of $K$.  Then for all $P$, $Q \in E(K^{\sep})$,
$$
        \deg (\pair{P}{Q}_\infty) = \pair{P}{Q}_{\NT},
$$
where here $\deg = -\ord_\infty$.
\end{thm}

\begin{proof}
Because the value of~$\Hhat_\infty$ is independent of a
chosen Weierstrass equation, we are free to fix a Weierstrass equation
which is minimal at all places of~$K$ extending~$\infty$.
Furthermore, it suffices that the theorem holds on~$E_\infty(K)$.
Given the ways we have normalized our pairings, we need to show that
\begin{equation} \label{eqn:HP}
        \deg (\Hhat_\infty(P)) = 2\,\hhat_{\NT} (P).
\end{equation}
By the definition of~$\Hhat_\infty$ along with~(\ref{eqn:ordiP})
and~(\ref{eqn:rhoKv}), we have
$$
        \Hhat_\infty(P)^{[K:k]} = \prod_{w \mid \infty} \pos{
                \NN^{K_w}_{k_\infty} \sigma_{E/\!K_w,\omega}(P)^{-2}}
        \,\prod_{v \neq \infty}
        \prod_{w \mid v}
                \pi_v^{[\FF_w:\FF_v] \lambda_w(P)},
$$
and thus the terms from~(\ref{eqn:HP}) corresponding to the finite
places exactly match those from~(\ref{eqn:NTloc}).  Finally we need to
show for each $w \mid \infty$ that, before taking norms, we have
$\ord_w (\sigma_{E/\!K_w,\omega}(P)^2) = \lambda_w(P)$, which follows
{from}~(\ref{thm:MT}) and the ultrametric inequality.
\end{proof}

\section{Non-degeneracy of $v$-adic heights} \label{sec:nondeg}

A symmetric bilinear pairing of abelian groups $E \times E \to G$ with
$G$ uniquely divisible is said to be \emph{non-degenerate} if the
kernel consists only of the torsion elements of $E$.
Theorem~\ref{thm:NText} assures us that the $\infty$-adic height
pairing is non-degenerate on $E(K)$.

On the other hand, the values of the $v$-adic pairing (for~$v$ a
finite place of~$k$) are \mbox{$1$-units} in~$\CC_v$, and so \emph{a
priori} we can say little about the non-degeneracy of these heights.
In this section and the next we investigate conditions for the
non-degeneracy of the pairing $\pairing_v$ on $E(K)$, and we find
that, for~$E/k$ defined over the rational function field,
non-degeneracy on $E(K)$ can be proven in many general cases.
Additionally, in the next section we show that the induced pairing
$$
        \pairing_v \colon (E(K) \otimes \ZZ_p) \times
                        (E(K) \otimes \ZZ_p) \to \sbstar{\CC}{v},
$$
is also non-degenerate under these same hypotheses.

\begin{rmk}
By contrast, for analytic $p$-adic height pairings for elliptic curves
defined over number fields~\cite{maztat83,schn82}, it is a difficult
problem to determine non-degeneracy (see~\cite{schn85}), and it is not
known in general when such pairings are non-degenerate.
\end{rmk}

For this section and the next we restrict ourselves to the following
situation: $E/k$ is an elliptic curve defined over the rational
function field with ordinary reduction at both~$v$ and~$\infty$.  This
reduction assumption is quite weak, since if $E$ is not supersingular,
then all but finitely many places will have ordinary reduction.  We
will further fix a Weierstrass equation~(\ref{eqn:wefix}) which is
minimal at both~$v$ and~$\infty$.

Let $K/k$ be a global function field.  If $v$ (resp.~$\infty$)
ramifies in~$K$, then we will also assume that~$E$ has good reduction
at~$v$ (resp.~$\infty$).  This assumption ensures that our chosen
equation is minimal at all places above~$v$ and~$\infty$ in~$K$.

For a point $P \in E_v(K) \cap E_\infty(K)$ and a fixed positive
integer~$N$, $V_N\colon E^{(p^N)} \to E$ is the Verschiebung (the dual
of the $p^N$-th power Frobenius morphism), and we define
$$
        V_N^{-1}(P) \subset E^{(p^N)}(k^{\sep})
$$
to be the inverse image of the point~$P$ under~$V_N$, which consists
of~$p^N$ distinct points.

Let $w$ denote a place of~$K$ above either~$v$ or~$\infty$.  If $P \in
E_v(K) \cap E_\infty(K)$, then because the formal group
$\Ehat^{(p^N)}/\OO_w$ has height~$1$, there is a unique point,
\begin{equation} \label{eqn:Qw}
        Q_w = Q_w(P) = Q_{w,N}(P) \in V_N^{-1}(P),
\end{equation}
such that $Q_w \in \Ehat^{(p^N)}(\OO_w)$.  For different places~$w$
and~$w'$ (either both extending~$v$ or both extending~$\infty$), the
points~$Q_w$ and $Q_{w'}$ are related.  Indeed, if we fix a separable
closure~$k_v^{\sep}$ so that
$$
        k_v \subset K_w \subset k_v^{\sep},
$$
we see that the absolute value of~$w'$ on~$K$ is obtained through an
embedding $\tau\colon K \hookrightarrow k_v^{\sep}$, which induces an
embedding $\tau\colon K_{w'} \hookrightarrow k_v^{\sep}$.  For such an
embedding $\tau$, we then have $\tau Q_{w'}(P) = Q_w(\tau P)$.

We now introduce a hypothesis which will play a crucial role in the
proof of Theorem~\ref{thm:Hvneq1}, the main theorem of this section.
In Section~\ref{sec:kummer} we will investigate the generality in
which this hypothesis holds.

\begin{hyp} \label{hyp:PN}
Let $P \in E_v(K) \cap E_\infty(K)$ and $N \geq 1$.  For fixed $w \mid
v$ and $\infbar \mid \infty$, there exists
$\gamma \in \Gal(k^{\sep}/k)$ so that for all embeddings $\tau \colon
K \hookrightarrow k^{\sep}$, we have $\gamma Q_w(\tau P) =
Q_{\infbar}(\tau P)$.
\end{hyp}

\begin{thm} \label{thm:Hvneq1}
Let $E/k$ be an elliptic curve with ordinary reduction at~$v$
and~$\infty$.  Let $K/k$ be a global function field.  If~$v$
(resp.~$\infty$) is ramified in~$K$, we will further assume that~$E$
has good reduction at~$v$ (resp.~$\infty$).  Let $P \in E_v(K) \cap
E_\infty(K)$.  Suppose that for any~$N$ such that
$$
        p^N \nmid 2\,[K:k] \hhat_{\NT}(P),
$$
the point $P$ satisfies Hypothesis~\ref{hyp:PN}.  Then
$$
        \Hhat_v(P) = \pair{P}{P}_v \neq 1.
$$
\end{thm}

We begin with a series of lemmas needed for the proof of
Theorem~\ref{thm:Hvneq1}.  Our eventual method will be to compare the
values of~$\Hhat_v$ and~$\Hhat_\infty$ and to use the known
non-triviality of~$\Hhat_\infty$ to prove that $\Hhat_v$ is also
non-trivial.  This first lemma exhibits one of the peculiarities of
characteristic~$p$, especially when compared to number fields.  Its
proof is fairly standard (see~{\S}8.2 of~\cite{goss:FF}).

\begin{lem} \label{lem:pNth}
For any place~$v$ of~$k$, let $L/k$ be an algebraic extension such
that $L \subset k_v$.  Then for every integer~$N\geq 1$ the map
$$
        \frac{L^\times}{(L^\times)^{p^N}} \to
        \frac{\sbstar{k}{v}}{(\sbstar{k}{v})^{p^N}}
$$
is an injection.
\end{lem}

We take $\omega = dx/(2y + a_1x + a_3)$, $\alpha$ the Hasse invariant
of $E$, and $z = -x/y$ a fixed formal parameter for $\Ehat$.

\begin{lem} \label{lem:sigdefR}
Let $R$ be the ring $\FF_q[a_1, a_2, a_3, a_4, a_6][1/\alpha] \subset
k$.  Then for all places $w \mid v$ and $w \mid \infty$ of $K$, the
following statements hold.
\begin{alphenumerate}
\item The formal group law for $\Ehat/\OO_w$ with parameter~$z$ is
defined over~$R$.
\item For every $m \in \ZZ$, the division polynomial~$f_m$ is an
element of $R[x,y]$.
\item The sigma function $\sigma(z) = \sigma_{E/K_w,\omega}(z)$ is
given by the formula in~(\ref{thm:sigprod}), and moreover $\sigma(z)
\in \power{R}{z}$.
\end{alphenumerate}
\end{lem}

\begin{proof}
Recalling the ramification hypotheses on~$v$ and~$\infty$ set earlier
in the section, the equation for~$E$ will be minimal at all places
of~$K$ extending them, explaining~(a) and~(b).  Part~(c) is then a
restatement of~(\ref{thm:sigprod}).
\end{proof}

\begin{lem} \label{lem:sigcong}
For every $P \in \Ehat(\OO_w)$ and integer $N \geq 1$,
$$
        \sigma_{E/K_w,\omega}(P) \equiv f_{p^N} ( \mult{p^N}^{-1}
        (z(P))) \mod{(\sbstar{K}{w})^{p^N}}.
$$
\end{lem}

\medskip
\begin{proof}
Using~(\ref{thm:sigprod}), we find
$$
        \sigma_{E/K_w,\omega}(P) \equiv
                \prod_{n=1}^N f_p ( \mult{p^n}^{-1} (z(P))
                )^{p^{2n-2}} \mod{ (\sbstar{K}{w})^{p^N}},
$$
and then the result follows by induction on~(\ref{eqn:fpz}).
\end{proof}

For an alternate interpretation of $f_{p^N} ( \mult{p^N}^{-1}
(z(P)))$, we consider that, as a function of $z$, $f_{p^N}(z) =
g_{p^N} (z^{p^N})$ is a Laurent series in $z^{p^N}$.  In fact,
$g_{p^N} \in k(E^{(p^N)})$, and thus for $Q_w$ as in~(\ref{eqn:Qw}),
\begin{equation} \label{eqn:gpN}
        f_{p^N} ( \mult{p^N}^{-1} (z(P)))
                = g_{p^N} (V_N^{-1}(z(P)) ) = g_{p^N}(Q_w).
\end{equation}

\begin{proof}[Proof of Theorem~\ref{thm:Hvneq1}]
Taking $R = \FF_q[a_1,a_2,a_3,a_4,a_6][1/\alpha]$, we note that $R
\subset \OO_w$ for all $w \mid v$ and $w\mid \infty$.
By~(\ref{lem:sigdefR}) we see that
$$
        \sigma_{E/K_w,\omega}(z) = \sigma(z),
$$
where $\sigma(z)$ is given by the product in~(\ref{thm:sigprod}).

By the construction of $\Hhat_\infty$ we know that
$\Hhat_\infty(P)^{[K:k]} \in \sbstar{k}{\infty}$.  The fact that $p^N$
does not divide $2\,[K:k]\hhat_{\NT}(P)$ combined
with~(\ref{thm:NText}) then implies that
\begin{equation} \label{eqn:Hinfn1}
        \Hhat_\infty(P)^{[K:k]} \not\equiv 1
                \mod{(\sbstar{k}{\infty})^{p^N}}.
\end{equation}
Letting $\pos{\ }_v$ and $\pos{\ }_\infty$ denote the positive parts
with respect to~$v$ and~$\infty$, we then have
$$
        \Hhat_\infty(P)^{[K:k]} = H^0(P) \, \prod_{w\mid
                \infty} \pos{ \NN^{K_w}_{k_\infty}
                \sigma(z(P))^{-2}}_\infty \in \sbstar{k}{\infty},
$$
where $H^0(P)$ is the product of all the local factors at the finite
places.  Furthermore, from the construction of~$\Hhat_v$ we have that
$$
        \Hhat_v(P)^{[K:k]} = \pos{H^0(P)}_v \, \prod_{w\mid v}
                \pos{\NN^{K_w}_{k_v} \sigma(z(P))^{-2}}_v
                \in \sbstar{k}{v},
$$
where $H^0(P)$ is the \emph{same} in both equations.  As we will only
be interested in these quantities up to $p^N$-th powers, it is
important to note that, for all $x \in \sbstar{k}{v}$, $x \equiv
\pos{x}_v$ modulo $(\sbstar{k}{v})^{p^N}$, and likewise for
$\sbstar{k}{\infty}$.

Now let $M \subset \Gal(k^{\sep}\!/k)$ be a set of representatives of
the embeddings of $K \hookrightarrow k^{\sep}$.  Up
to~$(\sbstar{k}{v})^{p^N}$, we know from~(\ref{lem:sigcong}) that
\begin{align*}
        \prod_{w \mid v} \NN^{K_w}_{k_v} \sigma(z(P))
                & \equiv \prod_{w \mid v} \NN^{K_w}_{k_v} f_{p^N}(
                \mult{p^N}^{-1}(z(P)))
                        \mod{(\sbstar{k}{v})^{p^N}} \\
                & \equiv \prod_{\tau \in M} g_{p^N}( Q_w(\tau P) )
                        \mod{(\sbstar{k}{v})^{p^N}},
\end{align*}
where for some fixed place $w \mid v$ the point $Q_w(\tau P) \in
V_N^{-1}(\tau P)$ is defined in~(\ref{eqn:Qw}) and~(\ref{eqn:gpN}).
Likewise, we can show
$$
        \prod_{w \mid \infty} \NN^{K_w}_{k_\infty} \sigma(z(P))
                \equiv \prod_{\tau \in M} g_{p^N}( Q_{\infbar}(\tau P))
                        \mod{(\sbstar{k}{\infty})^{p^N}},
$$
for some fixed place $\infbar \mid \infty$ and corresponding
$Q_{\infbar}(\tau P) \in V_N^{-1}(\tau P)$.  By (\ref{eqn:Hinfn1}) we
then have
$$
        \Hhat_\infty(P)^{[K:k]} \equiv H^0(P) \, \prod_{\tau \in M}
                g_{p^N}(Q_{\infbar}(\tau P))^{-2} \not\equiv 1
                \mod{(\sbstar{k}{\infty})^{p^N}}.
$$
Likewise, we have
$$
        \Hhat_v(P)^{[K:k]} \equiv H^0(P) \, \prod_{\tau \in M}
        g_{p^N}( Q_w(\tau P) )^{-2} \mod{(\sbstar{k}{v})^{p^N}}.
$$
For both of these two congruences, the right-hand sides are actually
elements of $k^{\sep}$.  Now by Hypothesis~\ref{hyp:PN} the points
$Q_w(\tau P)$ and $Q_{\infbar}(\tau P)$ are simultaneously
$\Gal(k^{\sep}\!/k)$-conjugates, and since $\Hhat_\infty(P)^{[K:k]}$
is not a $p^N$-th power, we conclude that $\Hhat_v(P)^{[K:k]}$ is also
not a $p^N$-th power by~(\ref{lem:pNth}).  In particular,
$$
        \Hhat_v(P)^{[K:k]} \neq 1,
$$
and the theorem follows.
\end{proof}

\section{Conditions for non-degeneracy} \label{sec:kummer}

The main question of this section is the generality in which
Hypothesis~\ref{hyp:PN} holds for points in $E_v(K) \cap E_\infty(K)$.
As in the previous section, we restrict ourselves to elliptic
curves~$E$ defined over the rational function field~$k$ with ordinary
reduction at both~$v$ and~$\infty$.  We will assume that $K/k$ is
Galois.  We let $\sT_p(E) = \varprojlim\, \ker V_n$ be the $p$-adic
Tate module.  Since~$E$ is ordinary, we have~$\sT_p(E) \cong \ZZ_p$ as
groups.  Our main result and an immediate corollary are as follows.

\begin{thm} \label{thm:Hvnev1}
  Let $E/k$ be an elliptic curve with ordinary reduction at~$v$
  and~$\infty$.  Let~$K/k$ be a finite Galois extension.  If~$K$ is
  ramified at~$v$ (resp.~$\infty$), then we further assume that~$E$
  has good reduction at~$v$ (resp.~$\infty$).  Finally, we assume that
  $\Gal(K^{\sep}\!/K) \to \Aut(\sT_p(E))$ is surjective and that for
  all places~$w$ extending~$v$ and~$\infty$ the reduction is good or
  nonsplit and the $\FF_w$-rational $p$-torsion on the reduction
  at~$w$ is trivial.  Then
\begin{alphenumerate}
\item $\Hhat_v(P) \neq 1$ for all non-torsion $P \in E(K)$.
\item For all $P$, $Q \in E(K)$, if $\pair{P}{Q}_{\NT} \neq 0$, then
$\pair{P}{Q}_v \neq 1$.
\end{alphenumerate}
\end{thm}

\begin{cor} \label{cor:pairvnond}
Under the conditions of Theorem~\ref{thm:Hvnev1}, the canonical height
pairing
$$
        \pairing_v\colon E(K) \times E(K) \to \sbstar{\CC}{v}
$$
is non-degenerate.
\end{cor}

The hypotheses in the above theorems are fairly weak, and thus the
results hold in some generality.  For example, the surjection of the
Galois representation is akin to the classical result that an integer
is a primitive root modulo~$p^n$, $n \geq 2$, if and only if it is a
primitive root modulo~$p^2$.  Similarly, it is easy to show that
$\Gal(K^{\sep}\!/K) \to \Aut(\sT_p(E))$ is surjective if and only if
it surjects onto $\Aut(\ker V_2)$.

Let $M \subset E(K)$ be a free abelian group and consider the tower of
field extensions
$$
        K \subset K_N \subset L_N,
$$
where $K_N = K(\ker V_N)$ and $L_N = K_N(V_N^{-1}(M))$.  If $M =
\sum_{i=1}^r \ZZ P_i$ has rank~$r$, there is a natural map
\begin{equation} \label{eqn:gal}
        \Gal(L_N/K_N) \to (\ker V_N)^r,
\end{equation}
via the Kummer pairing.  That is, $\tau \mapsto (\tau Q_i - Q_i)$,
where $Q_i \in V_N^{-1}(P_i)$.

The following proposition is a Verschiebung analogue of
Bash\-ma\-kov's Theorem~\cite{bash70,bash72}, and it uses a
Verschiebung descent to show that under certain conditions the image
of~(\ref{eqn:gal}) is as large as possible.  Since the Verschiebung is
separable, such descents behave much like prime-to-$p$ descents, as
opposed to the more delicate full $p$-descent in characteristic~$p$
(see~\cite{broumas97,ulmer91,voloch90}).  Moreover, the proof is
virtually identical to the one found in~{\S}V.5 of~\cite{lang:ECDA},
or see~\cite{map:PhD} for further details.

\begin{prop} \label{prop:bash}
Suppose $\Gal(K^{\sep}\!/K) \to \Aut(\sT_p(E))$ is surjective, and let
$M \subset E(K)$ be a torsion free subgroup of rank $s$.
\begin{alphenumerate}
\item Let $M = \ZZ P$.  If $P \notin V(E^{(p)}(K))$, then
        $\Gal(L_N/K_N) \cong \ker V_N$.
\item If $M \cap V(E^{(p)}(K)) = pM$,
        then $\Gal(L_N/K_N) \cong (\ker V_N)^s$.
\item In general, the Kummer pairing induces an isomorphism
$$
  \Gal(L_N/K_N) \cong \mathrm{Hom} \left(
    \frac{M}{M \cap V_N(E^{(p^N)}(K))}, \ker V_N \right).
$$
\end{alphenumerate}
\end{prop}

\begin{proof}[Proof of Theorem~\ref{thm:Hvnev1}]
It suffices to prove the theorem for points~$P$, $Q$ in $E_v(K) \cap
E_\infty(K)$.  For~(a), choosing an~$N$ so that
$$
        p^N \nmid 2\,[K:k]\hhat_{\NT}(P),
$$
we need to show that $P$ satisfies Hypothesis~\ref{hyp:PN}.  Fixing
places $w \mid v$ and $\infbar \mid \infty$, we need to show that for
all embeddings $\tau \colon K \hookrightarrow k^{\sep}$, the points
$Q_w(\tau P)$ and $Q_{\infbar}(\tau P)$ are simultaneously
$\Gal(k^{\sep}\!/k)$-conjugate.

Let $M = \sum \ZZ \tau P$, which is torsion-free by the hypotheses on
$P$.  By the elementary divisors theorem there is a $\ZZ$-basis
$\{R_i\}$ for $M$ so that
$$
  \frac{M}{M \cap V_N(E^{(p^N)}(K))} \cong \bigoplus
    \frac{\ZZ R_i}{p^{N-k_i}\ZZ R_i}
    \quad \text{for $0 \leq k_i \leq N$}.
$$
Note that as $p^{N-k_i} R_i \in M \cap V_N(E^{(p^N)}(K))$, our
hypothesis that the $p^k$-torsion $E^{(p^k)}(K)[p^k] = \{ O \}$ for
all $k$ implies that $R_i \in V_{k_i}(E^{(p^{k_i})}(K))$.  Thus via
the Kummer pairing with chosen basis $M = \sum \ZZ R_i$,
(\ref{prop:bash}c) becomes
\begin{equation} \label{eqn:Gker}
  \Gal(L_N/K_N) \cong \prod \ker V_{N,k_i},
\end{equation}
where $V_{N,k}$ is the kernel of the Verschiebung $E^{(p^N)} \to
E^{(p^k)}$.

By expressing each $\tau P$ as a linear combination of the $R_i$, we
need to show that there is a $\gamma \in \Gal(k^{\sep}\!/k)$ so that
\begin{equation} \label{eqn:simRi}
        \gamma Q_w(R_i) = Q_{\infbar}(R_i), \quad \text{for all $i$.}
\end{equation}

First, if $M \cap V(E^{(p)}(K)) = pM$, we apply (\ref{prop:bash}b).
Then there is an element $\gamma \in \Gal(k^{\sep}\!/K_N)$
accomodating Hypothesis~\ref{hyp:PN}, which we lift back to
$\Gal(k^{\sep}\!/k)$.

Second, if $M \cap V_N(E^{(p^N)}(K)) = M$, then for each $i$ we set
$Q(R_i)$ to be the unique element of $E^{(p^N)}(K) \cap
V_N^{-1}(R_i)$.  Indeed, because the $p^N$-torsion $E^{(p^N)}(K)[p^N]
= \{ O\}$ by assumption, these two sets meet at only one point.  Using
the assumptions on the $p$-torsion on the reductions, we have that
$Q(R_i) = Q_w(R_i) = Q_{\infbar}(R_i)$ for each $i$.

Finally, in the general case we proceed by combining the above two
arguments.  For each $i$, $k_i$ is a largest integer $k$ for which
$R_i \in V_k(E^{(p^k)}(K))$.  As in the preceding paragraph we have
$V_{N,k_i}(Q_w(R_i)) = V_{N,k_i}(Q_{\infbar}(R_i))$.  Thus from
(\ref{eqn:Gker}) we find that there is then an automorphism $\gamma
\in \Gal(k^{\sep}\!/k)$ which satisfies (\ref{eqn:simRi}).

For (b), the proof is much the same as in part (a).  By taking $N$ so
that $p^N \nmid 2\,[K:k] \pair{P}{Q}_{\NT}$, an argument similar to
the proof of~(\ref{thm:Hvneq1}) shows that we need to have $P$ and $Q$
satisfy Hypothesis~\ref{hyp:PN} simultaneously.  That is, we require
that there exist $\gamma \in \Gal(k^{\sep}\!/k)$ so that
$$
        \gamma Q_w(\tau P) = Q_{\infbar}(\tau P) \quad \text{and} \quad
        \gamma Q_w(\tau Q) = Q_{\infbar}(\tau Q)
$$
for all $\tau\colon K \hookrightarrow k^{\sep}$.  This is achieved by
applying the arguments of part (a) to $M = \sum \ZZ \tau P + \sum \ZZ
\tau Q$.
\end{proof}

\begin{rmk} \label{rmk:pow}
It should be noted that in~(\ref{thm:Hvneq1}) and~(\ref{thm:Hvnev1})
we have actually proven more than non-triviality.  We have shown,
under the hypotheses of these theorems, that for $P$, $Q \in E_v(K)
\cap E_\infty (K)$, if $p^N \nmid 2\,[K:k] \pair{P}{Q}_{\NT}$, then
$$
        \pair{P}{Q}_v^{[K:k]} \not\equiv 1 
        \mod{(\sbstar{k}{v})^{p^N}}.
$$
\end{rmk}

As pointed out by the referees, the non-degeneracy statement
in~(\ref{cor:pairvnond}) is not strong enough for most conceivable
applications.  Certainly, what would be better is a statement about
the non-triviality of a determinant for the pairing, a difficult
question to formulate since the values are taken in the multiplicative
group of $\CC_v$.  However, we prove the following corollary, which
shows that, under the hypotheses of~(\ref{thm:Hvnev1}), the pairing
$\pairing_v$ is non-degenerate on $E(K) \otimes \ZZ_p$.  This can be
used to show the non-existence of universal norms coming from a
certain pro-$p$-extension of $K$ (see~\cite{map:Univ,map:PhD}).

\begin{cor}
Under the conditions of Theorem~\ref{thm:Hvnev1}, the pairing
$$
        \pairing_v\colon (E(K) \otimes \ZZ_p) \times
        (E(K) \otimes \ZZ_p) \to \sbstar{\CC}{v}
$$
is non-degenerate.
\end{cor}

\begin{proof}
Because the determinant of the N\'eron-Tate height is non-zero, we
know $\pairing_{\NT}$ is non-degenerate on $E(K) \otimes \ZZ_p$.  Thus
we suppose $\pair{P}{Q}_{\NT} \neq 0$ for some $P$, $Q \in E(K)
\otimes \ZZ_p$.  By taking suitable multiples of $P$ and $Q$ we can
assume that both are in \mbox{$(E_v(K) \cap E_\infty(K)) \otimes
\ZZ_p$,} and then $p^N \nmid 2\,[K:k]\pair{P}{Q}_{\NT}$ for all $N$
large enough.  We write
$$
        P = P_1 + p^N P_2 \quad \text{and} \quad Q = Q_1 + p^N Q_2,
$$
where $P_1$, $Q_1 \in E_v(K) \cap E_\infty(K)$.  Then $p^N
\nmid 2\,[K:k]\pair{P_1}{Q_1}_{\NT}$, and by~(\ref{rmk:pow}), we have
that $\pair{P_1}{Q_1}_v^{[K:k]}$ is not a $p^N$-th power
in~$\sbstar{k}{v}$.  Moreover,
$$
        \pair{P}{Q}_v = \pair{P_1}{Q_1}_v \pair{P_2}{Q_1}_v^{p^N}
                \pair{P_1}{Q_2}_v^{p^N}\pair{P_2}{Q_2}_v^{p^{2N}},
$$
and we are done.
\end{proof}

\section{An example} \label{sec:exmp}

The Mazur-Tate sigma function is effectively computable
using~(\ref{thm:sigprod}).  It is thus possible to compute values
of~$\Hhat_v$ for a given elliptic curve and rational point.  For
example, suppose $E/k$ is defined over the rational function field and
has ordinary reduction at a finite place~$v$.  Then on a Weierstrass
equation minimal at all finite places, the height of a point $P \in
E_v(k)$ is given by
$$
        \Hhat_v(P) = \pos{\frac{\den (x(P))}{\sigma_{E/k_v,\omega}
                (P)^2} }_v,
$$
where $\den(x(P))$ is the denominator of the $x$-coordinate of~$P$.
Thus the $v$-adic height of $P$ can be obtained by computing values
of~$\sigma_{E/k_v,\omega}$.

Consider the elliptic curve $E/\FF_3(t)$ given by the equation
$$
        y^2 = x^3 + (t^2 - 1)x^2 + (t-1)^2(t^2-t-1)^2,
$$
which has discriminant $\Delta = (t+1)^3(t-1)^5(t^2-t-1)^2$ and Hasse
invariant $\alpha = t^2-1$.  In particular, $E$ has ordinary reduction
at~$t$.  By inspection~$E$ has two (linearly independent) $k$-rational
points
$$
        P = (0,(t-1)(t^2-t-1)), \quad
        Q = (t^2 - t, t^2 -1 ),
$$
neither of which is an element of~$E_t(k)$.  However calculations will
show that $30P$, $30Q \in E_t(k)$.  For $z = -x/y$ and $\omega =
dx/2y$, the first few terms of~$\sigma_{E,\omega}(z)$ are
\begin{eqnarray*}
        \sigma_{E,\omega}(z) & = &  z - (t^2-1)z^3 +
                \frac{(t-1)^2(t^2-t-1)^2}{t^2-1}\,z^5 \\
        & & {} - \frac{t^{11} + t^9 + t^5 - t^2 -1}{(t+1)^6}\,z^7
                +O(z^9)
\end{eqnarray*}
in $\laurent{\FF_3(t)}{z}$.  Accordingly,
$$
        \Hhat_t (30P) = 1 - t^3 - t^9 + t^{12} + t^{27} - t^{30} -
        t^{45} + O(t^{48})
$$
and
$$
        \Hhat_t(30Q) =  1 - t^3 + t^{18} - t^{21} + t^{27} +
        O(t^{30}).
$$

The curve~$E$ also has ordinary reduction at~$\infty$, and we can
calculate~$\Hhat_\infty$ similarly.  We find that $6P$, $6Q \in
E_\infty(k)$, and (with $t^{-1}$ the chosen uniformizer at~$\infty$)
$$
        \Hhat_\infty(6P) = t^{12} - t^{11} + t^{10} + t^8 + t^7 + t^6 + t^5
        + O(t^4)
$$
and
$$
        \Hhat_\infty(6Q) = t^{30} - t^{27} - t^{21} + t^{18} - t^3 + 1
        +O(t^{-3}).
$$
In terms of Hypothesis~\ref{hyp:PN}, it is worth noting that~$30Q$
satisfies the hypothesis, whereas $30P$ does not, although in both
cases the~$t$-adic height is non-trivial.

\section*{Acknowledgements}

This paper includes some of the results of the author's Brown
University Ph.D.\ thesis, and the author would especially like to
thank his advisor, Joseph~Silverman, for all of his encouragement and
guidance.  The author further thanks Antonios Broumas, Michael Rosen,
Jeremy Teitelbaum and Siman Wong for many helpful discussions on the
contents of this paper.  Finally, the author thanks the referees for
making many useful comments and suggestions.

\bibliographystyle{amsplain}

\begin{thebibliography}{10}

\bibitem{bash70}
M.~Bashmakov, \emph{Un th\'eor\'eme de finitude sur la cohomologie des courbes
  elliptiques}, C.~R.~Acad.\ Sci.\ Paris S\'er.~A-B \textbf{270} (1970),
  A999--A1001.

\bibitem{bash72}
M.~Bashmakov, \emph{The cohomology of abelian varieties over a number field},
  Russian Math.\ Surveys \textbf{27} (1972), 25--70.

\bibitem{bern81}
D.~Bernardi, \emph{Hauteurs $p$-adiques sur les courbes elliptiques},
  S\'eminaire de Th\'eorie des Nombres, Paris 1979--1980 (M.-J. Bertin, ed.),
  Birkh\"auser, Boston, 1981, pp.~1--14.

\bibitem{broumas97}
A.~Broumas, \emph{Effective $p$-descent}, Compositio Math. \textbf{107} (1997),
  125--141.

\bibitem{goss:FF}
D.~Goss, \emph{Basic structures of function field arithmetic}, Springer-Verlag,
  Berlin, 1996.

\bibitem{katmaz}
N.~M. Katz and B.~Mazur, \emph{Arithmetic moduli of elliptic curves}, Princeton
  University Press, Princeton, 1985.

\bibitem{lang:ECDA}
S.~Lang, \emph{Elliptic curves: Diophantine analysis}, Springer-Verlag, Berlin,
  1978.

\bibitem{maztat83}
B.~Mazur and J.~Tate, \emph{Canonical height pairings via biextensions},
  Arithmetic and Geometry: Papers Dedicated to I.~R.~Shafarevich on the
  Occasion of His Sixtieth Birthday (M.~Artin and J.~Tate, eds.), vol.~I,
  Birkh\"auser, Boston, 1983, pp.~195--237.

\bibitem{maztat91}
B.~Mazur and J.~Tate, \emph{The $p$-adic sigma function}, Duke
  Math.~J. \textbf{62} (1991), 663--688.

\bibitem{MTT86}
B.~Mazur, J.~Tate, and J.~Teitelbaum, \emph{On $p$-adic analogues of the
  conjectures of {B}irch and {S}winnerton-{D}yer}, Invent.\ Math. \textbf{84}
  (1986), 1--48.

\bibitem{map:Univ}
M.~A. Papanikolas, \emph{Universal norms on abelian varieties over global
  function fields}, (in preparation).

\bibitem{map:PhD}
M.~A. Papanikolas, \emph{Canonical heights in characteristic~$p$},
  Ph.D. thesis, Brown University, Providence, Rhode Island, 1998.

\bibitem{pr82}
B.~Perrin-Riou, \emph{Descente infinie et hauteur $p$-adique sur les courbes
  elliptiques \`a multiplication complexe}, Invent.\ Math. \textbf{70}
  (1982/83), 369--398.

\bibitem{schn82}
P.~Schneider, \emph{$p$-adic height pairings {I}}, Invent.\ Math. \textbf{69}
  (1982), 401--409.

\bibitem{schn85}
P.~Schneider, \emph{$p$-adic height pairings {II}}, Invent.\ Math. \textbf{79}
  (1985), 329--374.

\bibitem{sil:ATAEC}
J.~H. Silverman, \emph{Advanced topics in the arithmetic of elliptic curves},
  Springer-Verlag, New York, 1994.

\bibitem{ulmer91}
D.~L. Ulmer, \emph{$p$-descent in characteristic $p$}, Duke Math.~J.
  \textbf{62} (1991), 237--265.

\bibitem{voloch90}
J.~F. Voloch, \emph{Explicit $p$-descent for elliptic curves in characteristic
  $p$}, Compositio Math. \textbf{74} (1990), 247--258.

\bibitem{voloch:zeta}
J.~F. Voloch, \emph{An analogue of the {W}eierstrass $\zeta$-function in
  characteristic~$p$}, Acta Arith. \textbf{79} (1997), no.~1, 1--6.

\end{thebibliography}

%\end{article}

\end{document}